\documentclass[review]{elsarticle}
\usepackage[bookmarks=true,colorlinks=true,urlcolor=blue,citecolor=red,linkcolor=blue,linktocpage,pdfpagelabels,bookmarksnumbered,bookmarksopen]{hyperref}
\usepackage{lineno,hyperref}

\modulolinenumbers[5]

\journal{arXiv}
\usepackage{amscd}
\usepackage{mathrsfs}
\usepackage{amsfonts}
\usepackage{graphicx}
\usepackage{amsmath,amscd,stmaryrd,amssymb,array, amsfonts,amsmath,amssymb}
\usepackage{enumitem}
\usepackage{color}
\usepackage{verbatim}
\usepackage{extarrows}

\numberwithin{figure}{section}
 \numberwithin{equation}{section}

\newtheorem{theorem}{Theorem}[section]
\newtheorem{proposition}[theorem]{Proposition}
\newtheorem{definition}[theorem]{Definition}

\newtheorem{lemma}[theorem]{Lemma}
\newtheorem{remark}[theorem]{Remark}

\newcommand{\bu}{{u}}
\newcommand{\bbf}{{f}}
\newcommand{\bg}{{g}}
\newcommand{\bh}{{h}}

\newcommand{\bV}{{ V}}

\newcommand{\bH}{{ H}}

\newcommand{\cA}{{\mathcal A}}
\newcommand{\cB}{{\mathcal B}}

\newcommand{\cC}{{\mathcal C}}
\newcommand{\cD}{{\mathcal D}}
\newcommand{\cF}{{\mathcal F}}

\newcommand{\cL}{{\mathcal L}}

\newcommand{\sB}{{\mathscr B}}

\newcommand{\sD}{{\mathscr D}}

\newcommand{\sS}{{\mathscr S}}

\newcommand{\sX}{{\mathscr X}}

\def\be{\begin{equation}}
\def\ee{\end{equation}}
\def\bes{\begin{equation*}}
\def\ees{\end{equation*}}
\def\bsp{\begin{split}}
\def\esp{\end{split}}

\def\ba{\begin{array}}
\def\ea{\end{array}}
\def\benu{\begin{enumerate}}
\def\eenu{\end{enumerate}}
\def\bt{\begin{theorem}}
\def\et{\end{theorem}}
\def\bp{\begin{proposition}}
\def\ep{\end{proposition}}
\def\bl{\begin{lemma}}
\def\el{\end{lemma}}
\def\br{\begin{remark}}
\def\er{\end{remark}}
\def\bd{\begin{definition}}
\def\ed{\end{definition}}

\def\b{\beta}
\def\De{\Delta}

\def\pa{\partial}

\def\lam{\lambda}
\def\Lam{\Lambda}

\def\sig{\sigma}

\def\gam{\gamma}

\def\a{\alpha}
\def\W{\Omega}

\def\.{\cdot}

\def\R{\mathbb{R}}

\def\A{\forall}
\def\ol{\overline}

\def\Cap{\bigcap}\def\Cup{\bigcup}

\def\ra{\rightarrow}

\def\~{\tilde}
\def\8{\infty}
\def\X{\times}

\def\mb{\mbox}

\def\llg{\left\langle}
\def\rrg{\right\rangle}

\def\Hs{\hspace{1cm}}\def\hs{\hspace{0.5cm}}
\def\Vs{\vskip8pt}\def\vs{\vskip4pt}

\def\({\left(}\def\){\right)}

\begin{document}

\begin{frontmatter}

\title{ Global Existence, Regularity, and Dissipativity of Reaction-diffusion Equations with State-dependent Delay and  Supercritical Nonlinearities\footnote{This work was supported by the National Natural Science Foundation of China [12271399, 11871368]. }
}


\author[mymainaddress]{Ruijing Wang}
\ead{wrj\_math@163.com}
\address[mymainaddress]{School of Science, Qingdao University of Technology, Qingdao 266520,  China}

\author[mysecondaryaddress]{Desheng Li\corref{mycorrespondingauthor}}
\cortext[mycorrespondingauthor]{Corresponding author}
\ead{lidsmath@tju.edu.cn}
\address[mysecondaryaddress]{School of Mathematics,  Tianjin University, Tianjin 300350,  China}

\begin{abstract}
This work aims to study  the initial-boundary  value problem of the reaction-diffusion  equation
$\pa_{t}u-\Delta u=f(u)+g(u(t-\tau(t,u_t)))+h(t,x)$ in a bounded domain with  state-dependent delay and supercritical nonlinearities. We  establish the global existence and discuss the regularity and  dissipativity of the problem under weaker assumptions. In particular,  the existence of a global pullback attractor is proved regardless of uniqueness.
\end{abstract}

\begin{keyword}
Reaction-diffusion  equation, state-dependent delay \sep supercritical nonlinearity\sep  global existence\sep regularity \sep pullback attractor.
 \MSC[2020]   35A01  \sep 35B40 \sep 35B41 \sep 35B65.
\end{keyword}

\end{frontmatter}


\section{Introduction}\label{s:1}

The study of functional differential equations with SDDs (state-dependent delays) can be traced back to the earlier work of Possion in 1806 \cite{Poi}. In 1963,  Driver \cite{Dri} introduced a  formula involving the charge field into the Lorentz force law. This gives the birth  to  the first intuitive  functional differential equation with SDD. Since then various  differential systems with SDDs were proposed in different  areas, which coursed an increasing interest in the investigation of qualitative properties  of such systems; see e.g.   \cite{Aie,Cal,Har,Her1,Her5,LXD} et al.

Comparing with constant or time-dependent delay  differential systems, generally the research of systems with SDDs  encounters more difficulties, due to the lack of  smoothness in the solution spaces. In spite of this severe fact, one can find a large number of beautiful  works addressing this challenging problem \cite{Hartung,Kri1,Kri,Lv,Rez3,Rez2,Wal2}.
For instance, in \cite{Kri1}, the authors studied
the delay ODE $\dot{x}(t)=-\mu x(t)+F(x(t-\tau(x(t))))$ and described the asymptotic behavior of the slowly oscillating solutions, which is based on the sharp condition $F'<0$.
Hartung \cite{Hartung} obtained the existence and uniqueness of a solution to the PDE $ \pa_t u=F(t,u_t,u(t-\tau(t,u_t)))$ and furthermore provided the differentiability of the solution with respect to the initial function, but only if the delay function $\tau$ is strictly increasing.
These are  interesting and elegant conditions, but also important restrictions on applicability.

This paper is devoted to the following non-autonomous  reaction-diffusion equation  in a smooth bounded domain $\W\subset \mathbb{R}^d$ with fast-growing nonlinearities:
\be\label{e:1.1}\left\{\ba{ll}
\pa_{t}u-\Delta u=f(u)+g(u(t-\tau(t,u_t)))+h(t,x),\\[1ex]
u|_{\pa\W}=0, \\[1ex]
u|_{\overline{\W}\times [-r,0]}=\phi,\ea\right.
\ee
where $f(\cdot), g(\cdot)\in C^1(\mathbb{R})$, $ 0\leq \tau(\cdot,\cdot)\leq r$, and $u_t$ denotes the history of the state $u(\cdot)$ at the time $t$ is given by
\be\label{e:1.5}
u_t(\theta)=u(t+\theta), \Hs \theta \in [-r,0].
\ee
We are basically interested in the case where $f$ and $g$ are supercritical with  the delay function $\tau$ satisfying some weaker assumptions than those proposed in the literature.

 System  \eqref{e:1.1} plays a   crucial role in the theory of functional PDEs. In many cases $f$ and $g$ are  assumed to have polynomial growth rates:
\benu
  \item[{\bf(H1)}] There exist $p,\beta, a_0, b_0>0$ such that
$$
|f(s)|\leq a_0(|s|^p +1),\Hs s\in \mathbb{R},
$$
$$
|g(v)|\leq b_0(|v|^\beta +1),\Hs v\in \mathbb{R}.
$$
\eenu
Generally, there exists a critical value $p_c>1$, based on which equation \eqref{e:1.1} can be divided into three distinct cases: the subcritical case ($\max\{p,\b\}<p_c$), the critical case ($\max\{p,\b\}=p_c$), and the supercritical case ($\max\{p,\b\}>p_c$).
In the  subcritical and critical cases one can find numerous works dealing with the local or global well-posedness of the system.
For instance, for the equation
$$
\pa_t u -\Delta u(t)+F(u(t-\tau(u_t)))+G(u(t))=h,\Hs t>0,
$$
if $F$ and $\tau$ are globally Lipschitz and $G$ is locally Lipschitz, one can find some nice results concerning the well-posedness of solutions and the existence of attractors in Chueshov et al.\cite{Chu}. Closely related works can be found in \cite{Her5,Her3,Kri,Lv,Travis2, Hale2,Hale1,Kuang,Tra,Wu} etc., where some  additional requirements on the nonlinear and delay terms are requited to ensure the well-posedness.
We mention that these works are mainly focused on the case where $F$ is subcritical with the delay term being at most sublinear.


In the supercritical case, the  discussions about the global   well-posedness  and asymptotic behaviour of the problem becomes quite   complicated, due to the lack of  appropriate  functional framework and the possible existence of finite-time blow up  of solutions.
 For the constant delay or time-dependent delay equations, dissipative type structure conditions have received much attention in recent years as natural candidates for solving the global well-posedness problems; see e.g. \cite{Kostianko,Li2,LIX,Zelik2000} and references therein. Inspired by these works, our results are established based on the following dissipative condition, which allows $f$ to have an arbitrary polynomial growth rate:
\benu
  \item[{\bf(H2)}] There exist $\a\geq 1$, $N \geq 0$ and $\Lam>0$ such that
$$
f(s)s\leq -\Lam |s|^{\a+1}+N,\Hs s\in\mathbb{R}.
$$
\eenu
According to this, we impose the growth requirement on $g$, which reflects the a competition balance between the two nonlinearities:
\benu
  \item[{\bf(H3)}] $\b < \a$, where $\b$ is the exponent of $g$ given in  (H1).
\eenu

The study of qualitative properties of differential equations with SDDs is much more difficult than those with constant or time-dependent delays. To assume well-posedness,    in \cite{Her4} Hernandez et al. assumed that  $\tau$ to satisfy $L^q$-Lipschitz continuity along with some other hypotheses. Here we use some weaker assumptions. Specifically, we only assume that  $\tau$ satisfies the assumption:
\benu
  \item[{\bf(H0)}]
If   $\varphi_k\rightharpoonup\varphi$  in $H_0^1(\W)$  as $k\ra\8$, then 
$$\tau(t,\varphi_k)\rightarrow\tau(t,\varphi),\Hs k\ra\8,$$
\eenu

Now let us give a brief description to our main results.
Define
\bes\label{e:x.1.1}
 p_0=\frac{(\a-1)}{\a-\b}\b,\hs q_0=\mb{max}\{2p,2\b,p_0\},
\ees
where $\a$, $\b$ and $p$ are positive constants mentioned in (H1)-(H3).  For  $q\geq1$, we also write  $q_{\a}:=q-1+\a$.

For a  Banach space $X$, we denote
$$\cC_X=C([-r,0];X),\hs \cL^\8_X=L^\8(-r,0;X),$$
where $r$ is the upper bound of the delay functions $\tau$ given  in \eqref{e:1.5}, and let
$$
\bH=L^2(\W),\hs \bV_1=H^1_0(\W),\hs \bV_2=H^2(\W)\cap H^1_0(\W),
$$
and
\bes
\sX_i^q=\cC_{\bV_i}\cap \cL^\8_{L^q(\W)},\Hs i=1,2.
\ees
Besides, denote $H^{\zeta}:=D((-\Delta)^{\zeta})=H_0^{\zeta}(\W)\cap H^{2\zeta}(\W)$ by the $\zeta$-fractional power of $H$ with the norm $\|x\|_{\zeta}=\|(-\Delta)^{\zeta} x\|$, where $0\leq\zeta\leq 1$.

Our work focuses on proving the global existence, regularity of solutions for \eqref{e:1.1} and the existence of global pullback attractors by the Galerkin method and the semigroup approaches, only if the delay function $\tau$ is weakly continuous (satisfies (H0)). These results are summarized in the following theorems.
\bt\label{t:x.1.1}
Assume $\tau$ satisfies (H0), $f$ and $g$ satisfy (H1)-(H3). Let $q_0<q<\8$, and $ h\in L^{\8}(\mathbb{R};L^{q_\a/\a}(\W))\cap C(\mathbb{R};L^{q_\a/\a}(\W))$. If $\phi\in \mathscr{X}_1^q$, then the solutions $u=u(t;\phi)$ of \eqref{e:1.1} satisfy that
\be\label{e:1.22}
\left\{\ba{ll}  u\in C\([-r,\8); V_1\)\cap L^\8\(-r,\8; V_1\)\cap L^\8\(-r,\8;L^q(\W)\),\\[1.5ex] u\in L^2\(0,T;V_2\),\hs  u'\in L^2\(0,T;H\),\Hs\A\,T>0,\\[.5ex]\ea\right.\ee
Furthermore, there exist $B_i$, $\eta_i$ and  $\rho_i>0$ ($i=0,1$) such that
\be\label{e:q.3}
|u (t;\phi )|_{q}\leq B_0 \|\phi\|_{\cL_{L^q(\W)}^\8} e^{-\eta_0 t} +\rho_0,\Hs t\geq 0, \hs 1<q<\8,
\ee
\be\label{e:q.5}
\| u (t;\phi )\|_1 \leq B_1\(\|\phi\|_{\cC_{V_1}}+\|\phi\|_{\cL_{L^q(\W)}^\8}^{q/2}\) e^{-\eta_1 t}+\rho_1,\Hs  \,t\geq 0.
\ee
\et

\begin{theorem}\label{t:1.2x} Assume that hypotheses {\em (H0)-(H3)} are fulfilled. Suppose
$\bh\in  L^{\8}(\mathbb{R};L^{\8}(\W))\cap C(\mathbb{R};L^{\8}(\W))$.
If $\phi\in \sX_1^\8$, then  {\em \eqref{e:1.1}} has solutions $\bu=\bu(t;\phi)$ satisfying \eqref{e:1.22}  with $q=\8$; furthermore, there exist $\rho_*,\nu_*>0$ such that
\be\label{e:q.4}
\ba{lll}
| u (t;\phi )|_{\8}&\leq
\left\{\ba{ll} \|\phi\|_{\cL_{L^\8(\W)}^\8}+\rho_*,\hs &\mb{if }\,r>0;\\[1ex]
\|\phi\|_{\cL_{L^\8(\W)}^\8}e^{-\nu_*t}+\rho_*,\hs &\mb{if }\,r=0,\\[1ex]
\ea\right.
\ea
\ee
\end{theorem}

Thanks to the following compact embedding, which was described in more detail by Henry in \cite[Chap.1]{Henry}:
\be\label{e:q.2}
H_0^1(\W)\cap H^2(\W) \hookrightarrow H^{\zeta_2} \hookrightarrow H_0^1(\W)\hookrightarrow H^{\zeta_1} \hookrightarrow L^2(\W),\Hs 0<\zeta_1<1/2<\zeta_2<1,
\ee
we prove the boundedness and equi-continuity  for the solutions of \eqref{e:1.1} in the fractional power space $H^{\zeta} (\frac{1}{2}<\zeta<1)$:
\bes
\begin{split}
\|u\|_{L^{\8}(r,\8;H^{\zeta})}&\leq M^* \|\phi\|_{\cL_{L^q(\W)}^\8} + C_b.\\[.5ex]
\|u(t+\nu)-u(t)\|_{H^1_0(\W)}&\leq L \, (\nu^{\zeta-\frac{1}{2}}+\nu^{\frac{1}{2}}),\hs \forall \nu>0,\,\, t\geq 2r+\eta.
\end{split}\ees
Obviously these important conclusions imply the asymptotic compactness of the set-valued solution semigroup. On the other hand, the concept of a pullback attractor is also suitable for describing its global dynamics for the nonautonomous system \eqref{e:1.1}. We refer the interested reader
to \cite{Carva,Che,Caraballo} etc. for general theories on this topic. Here we borrow some ideas from \cite{Car1} and \cite{Coti}, and carry out our work in the context combining set-valued dynamical systems with pullback $\sD$-attractors.

Consider the family of systems corresponding to \eqref{e:1.1}  with the time symbol  $\sigma\in \mathbb{R}$, whose abstract form is given by:
\be\label{e:xx2}
\left\{
  \begin{array}{ll}
    \pa_t u-\Delta u=f(u)+g(u(t-\tau(t+\sigma,u_t)))+h(t+\sigma, x), \Hs t>0,\\[1ex]
u_0=\phi.\\[.5ex]
  \end{array}
\right.
\ee
Clearly Theorem \ref{t:x.1.1} and Theorem \ref{t:1.2x} remain true for system \eqref{e:xx2} with the constants therein being independent of $\sigma$ (see Section \ref{s:3}-\ref{s:4} for details). Define
$$
\Phi=\Phi(t,\sigma,\phi)=\{\bu_t\,|\,\bu\in \sS(\phi)\},
$$
where $\sS(\phi)$  is the set of solutions of \eqref{e:xx2} with initial value $\phi$.
Let $\sD$ be the family of bounded sets $\sB^q(R)$ $(R>0)$ given by
$$
\sB^q(R)=\{\phi\in \sX_1^q \,|\, \|\phi\|_{ \mathcal{C}_{V_1}}+\|\phi\|_{\mathcal{L}^{\8}_{L^q (\W)}}\leq R\}.
$$
Then we have following results concerning the existence of pullback attractors:
\bt\label{t:1.1x}
In addition to {\em (H0)}-{\em (H3)},
assume $f$ and $g$ satisfy
 \benu
 \item[{\bf(H3)}] $  \bbf(0)=0=\bg(0)$.
 \eenu
Let $ h\in L^{\8}(\mathbb{R};L^{q_\a/\a}(\W))\cap C(\mathbb{R};L^{q_\a/\a}(\W))$ and $q_0<q\leq \8$. Then  $\Phi$  has  a global pullback $\sD$-attractor $\cA$ in $\sX_1^q$
which pullback attracts each element in $\sD$ in the topology of $\mathcal{C}_{V_1}$.
\et

This paper is organized as follows. In Section \ref{s:2} we recall several functional spaces and give definitions of weak and mild solutions for \eqref{e:1.1}.
Section \ref{s:3} is devoted to the proof of Theorems \ref{t:x.1.1} and \ref{t:1.2x}.
In Section \ref{s:4}, we obtain some good results on the boundedness and equi-continuity of solutions in fractional power spaces. In Section \ref{s:5}, we give precise statements and detailed proofs  about the existence of the global pullback attractors of  \eqref{e:1.1}.

\section{Preliminaries}\label{s:2}
\subsection{Functional  spaces}\label{s:2.2}
\noindent $\bullet$ {\bf Sobolev spaces and  their associated functional spaces.}
The spaces $W^{k,p}$ and $W^{k,p}_0$ are the usual. Denote by $|\cdot|_p$ the usual norm on the Lebesgue space $L^p (\W)$, and let $(\cdot,\cdot)$ be the scalar product on $H=L^2(\Omega)$.
Let $H^k(\W)=W^{k,2}(\W)$, and $H^k_0(\W)=W^{k,2}_0(\W)$. Set the Hilbert spaces
$$
V_1=H^1_0(\W),\hs V_2=H^2(\W)\cap H^1_0(\W)
$$
with  the following scalar products respectively:
$$
\llg u,v \rrg_1=(\nabla u,\nabla v) \, (\forall u,v\in V_1),\hs \llg u,v \rrg_2=(  \Delta u,  \Delta v) \, (\forall u,v\in V_2).
$$
Denote by $\|\cdot\|_k$ ($k=1,2$) the norm induced by the scalar product $\llg \cdot,\cdot \rrg_k$. It is trival to check that $\|\cdot\|_k$ is a norm on $V_k$ equivalent to the usual $H^k(\W)$-norm.

Recall that $r$  is the number that appears in \eqref{e:1.5}. Given a Banach space $X$, set
\be\label{e:1.6}
\mathcal{C}_{X}=C([-r,0];X),\hs \mathcal{L}^{\infty}_{X}=L^{\infty}(-r,0;X).
\ee
These two spaces are equipped with the norms:
$$
\|\cdot\|_{\mathcal{C}_{X}}:=\|\cdot\|_{C([-r,0];X)},\hs \|\cdot\|_{\mathcal{L}^{\infty}_{X}}:=\|\cdot\|_{L^{\infty}(-r,0;X)}.
$$
If $r$ is taken specifically to be $0$, we assign $\mathcal{C}_{X}=\mathcal{L}^{\infty}_{X}=X$. For notation convenience, let
$$\label{e:1.7}
\mathscr{X}^q_i=\mathcal{C}_{V_i}\cap \mathcal{L}^{\infty}_{L^q(\W)},\Hs i=1,2.
$$

\noindent $\bullet$ {\bf Fractional power spaces }
We rewrite the equation in \eqref{e:1.1} along with   the Dirichlet boundary condition $u|_{\pa \W}=0$ into a weak abstract form:
\be\label{e:1.11}
u'+Au=f(u)+g(u(t-\tau(t,u_t)))+h,
\ee
where $A : V_1 \rightarrow V^{'}_1$ is the operator defined as
\be\label{e:1.11-1}
\llg Au,v \rrg=\int_{\W} \nabla u \cdot \nabla v dx,\Hs \forall u,v\in V_1.
\ee
The initial condition in \eqref{e:1.1} can be written as
\be\label{e:1.12}
u(s)=\phi(s),\Hs s\in [-r,0].
\ee

In some work concerning  the sectorial operator theory, it is easy to see that there exists $\delta>0$  such that $\mb{Re}\,\sig(A)>\delta$. We use the notation $A^{\zeta}$($ 0\leq \zeta< 1$) for the $\zeta$-fractional power  $A^{\zeta}: D(A^{\zeta})\mapsto H $ of $A$ and the symbol $H^{\zeta}$ for the domain of $A^{\zeta}$ endowed with the norm $\|x\|_{\zeta}=\|A^{\zeta} x\|$ ($x\in H^{\zeta}$), where $\|\cdot\|$ is the norm of $H$.
\bl\label{l:1.3}
{\em\cite{Henry}} Assign $H^0=H$ for $\zeta=0$, then $H^{\zeta}$ is a Banach space in the norm $\|\cdot\|_{\zeta}$ for $0 \leq \zeta< 1$; If $0\leq \zeta_1 < \zeta_2< 1$, $H^{\zeta_2}$ is a dense subspace of $H^{\zeta_1}$ with compact continuous inclusion, and it follows directly that
\be\label{e:2.4}
H^2(\W)\cap H_0^1(\W)\hookrightarrow H^{\zeta_2} \hookrightarrow H_0^1(\W) \hookrightarrow H^{\zeta_1}\hookrightarrow H,\Hs 0<\zeta_1<1/2<\zeta_2<1.
\ee
Clearly all of these embeddings are compact.
\el

\bl\label{l:1.2}
{\em\cite{Henry}} Recall that $\mb{Re}\,\sig(A)>\delta>0$, then for any $\zeta \geq 0$ there is $C_{\zeta}<\8$ such that
\be\label{e:1.14}
\|A^{\zeta} e^{-t A}\|\leq C_{\zeta} t^{-\zeta} e^{-\delta t},\Hs t\geq 0.
\ee
\el

\subsection{Definitions of weak solutions and mild solutions}\label{s:2.1}
\begin{definition}\label{d:1.1}
Let $J\subset \mathbb{R}$ be an interval. A function $u=u(t)$ taking values in $V_1$ is called a weak solution of \eqref{e:1.11} on $J$, if for any compact interval $[a,b]\subset J$,
\be\label{e:2.5}
\begin{split}
u\in C([a-r,b];H)\cap &L^2(a-r,b;V_1)\cap L^{2\hat{p}}(a-r,b;L^{2 \hat{p}}(\W)),
\\[.5ex]
&u'\in L^2(a,b;V'_1),
\end{split}
\ee
where $\hat{p}=\max\{p,\b\}$.
Furthermore, the following equation holds in the distribution sense on $J$:
$$
\llg u',w\rrg+\llg Au,w\rrg=\int_{\W} (f(u)+g(u(t-\tau(t,u_t)))+h)\,w\, dx,\Hs \forall w\in V_1.
$$
\end{definition}
\begin{remark}\label{r:1.1}
If $u$ and $v$ are weak solutions of \eqref{e:1.11} on the intervals $[t_0, t_1]$ and $[t_1, t_2]$ respectively with $u(t)=v(t)$ for $t \in [t_1-r, t_1]$, then the union $z$ of $u$ and $v$ defined by
$$
z(t)=u(t)\,(t\leq t_1),\hs z(t)=v(t)\,(t>t_1)
$$
is a weak solution of \eqref{e:1.11} on $[t_0, t_2]$.
\end{remark}

Note that the requirements in \eqref{e:2.5} and the growth restrictions on the nonlinear terms in (H1) imply that $f(u),g(u(t-\tau(u_t)))\in L^2(a,b;H)$. Moreover, according to the standard theory of linear growth equations in Banach spaces (see \cite [Chap. 2, Section 3.2]{Tem}), it can be easily deduced that the weak solution $u$ of \eqref{e:1.11} on the interval $J$ must belong to $C(J,H)$.

\begin{definition}\label{d:1.2}
A weak solution $u$ of \eqref{e:1.11} on $(0,T)$ is called a weak solution of system \eqref{e:1.1}, if $u\in C([-r,T];H)$ and fulfills the initial condition \eqref{e:1.12}.
\end{definition}

\begin{definition}\label{d:1.3}
Let $T(t)$ be the $C_0$ semigroup generated by the   infinitesimal generator $-A$. We call a function $u\in C([-r,T];H)$  is the mild solution of \eqref{e:1.11} on $[-r, T]$, if $u$ is given by
\be\label{e:1.13}
u(t)=T(t)u(0)+\int^t_0 T(t-s) w(s)ds,\Hs 0\leq t\leq T,
\ee
where $w(t)=f(u(t))+g(u(t-\tau(t,u_t)))+h(t)$.
\end{definition}


\begin{theorem}\label{t:2.2}
{\em\cite{Pazy}} \eqref{e:1.13} holds if for any $u(0)\in L^2(\W)$, there exists a weak solution $u$ on the interval $[0, T]$ of the equation \eqref{e:1.11}, and $-A$ is an infinitesimal generator of $C_0$ semigroup $T(t)$.

\end{theorem}

From now on, we will simply call a \emph{weak solution} of \eqref{e:1.1} as a \emph{solution}, but the mild solution will not be abbreviated.

\section{Global $L^{q}$ and $H^1$-estimates of $u$}\label{s:3}
In this section we first derive some estimates for the Galerkin approximations $u_k$ of $u$,  and then employ compactness methods to obtain dissipative and regularity results of $u$ in $L^q(\W)$ and $H^1(\W)$.

\subsection{Global $L^{q}$ and $H^1$-estimates of $u_k$}
Let $\{w_j\}^{\8}_{j=1}$ be an orthogonal basis of $L^2$ consisting of eigenvectors of $A=-\Delta$. Given $\phi \in \mathcal{C}_{V_1}\cap \mathcal{L}^{\8}_{L^q (\W)}$, by the \cite[Theorem 7.1]{Li2} one can pick a sequence of smooth functions $\phi_k=\sum^k_{j=1}c_{kj}(t)w_j$\,($k=1,2,\cdots$) such that
\be\label{e:2.1}
\phi_k\rightarrow\phi \,\,( \mb{in}\,\, \mathcal{C}_{V_1}),\hs \|\phi_k\|_{\mathcal{L}^{\8}_{L^q (\W)}}\leq 8 \|\phi\|_{\mathcal{L}^{\8}_{L^q (\W)}}\,\, (k\in \mathbb{N}).
\ee
Similarly we can also take a sequence of functions $\{h_k\}$ such that $h_k \rightarrow h$ in appropriate spaces. For each $k$, let
\be\label{e:2.2}
u_k(t)=\sum^k_{j=1}a_{kj}(t)w_j
\ee
be a Galerkin approximation solution of \eqref{e:1.1},
\be\label{e:2.3}
\left\{
\begin{split}&\(\frac{d u_k}{dt}-\De u_k,w_j\)=\( f( u_k)+g(( u_k(t-\tau(t,(u_k)_t)))+ h_k,w_j\),\hs 1\leq j\leq k,\\
 & u_k(s)=\phi_k(s),\hs s\in[-r,0].
\end{split}
\right.
\ee
By the basic theory of the ODE, we know that $u_k$ is sufficiently regular so that all calculations can be performed rigorously on  $u_k$.




Define
\be\label{e:1.3}
q_0=\mb{max}\{2p,2\b,p_0\}, \Hs  \mb{where}\,\, p_0=\frac{(\a-1)}{\a-\b} \b.
\ee
Here $\a$, $\b$, and $p$ are the positive numbers involved   in (H1)-(H2). For $q\geq 1$, we also write
\bes\label{e:1.4}
q-1+\a:=q_{\a}.
\ees
It is standard to discuss the estimates of the Galerkin approximations $u_k$ using the Galerkin method; see e.g. \cite{Chu,Evans,Li2}.
Here we employ  similar computational techniques as described in \cite{Li} and \cite{Li2} to obtain corresponding estimates. Our results are summarized in the theorems below, with some proofs omitted.

\bl\label{t:1.1}
Assume that $\tau$ satisfies {\em(H0)}, $f$ and $g$ satisfy {\em (H1)-(H3)}. Let $q_0< q< \8$, and suppose that $  h\in L^{\8}(\mathbb{R};L^{q_\a/\a}(\W))\cap C(\mathbb{R};L^{q_\a/\a}(\W))$. Then there exist $M,\lam_q,R_q>0$ (where $M$ is independent of $q$), such that for all $\phi\in  \mathscr{X}_1^q$, the solution of \eqref{e:2.3} satisfies
\be\label{e:3.1}
|u_k(t;\phi_k)|_{q}^{q}\leq Me^{-\lam_q t}\|\phi\|_{\cL_{L^q(\W)}^\8}^{q}+R_q,\Hs t\geq -r,
\ee
where $\phi_k$ is the sequence shown in \eqref{e:2.1}.
\el

\bl\label{t:1.2}
Assume that $\tau$ satisfies {\em(H0)}, $f$ and $g$ satisfy {\em (H1)-(H3)}. Let $\bh\in  L^{\8}(\mathbb{R};L^{\8}(\W))\cap C(\mathbb{R};L^{\8}(\W))$. Then there exist $\rho_*,\lam_*>0$  such that for all $\phi\in \mathscr{X}_1^\8$, the solution of \eqref{e:2.3} satisfies
\be\label{e:3.2}\ba{lll}
| u_k(t;\phi_k)|_{\8}&\leq
\left\{\ba{ll} \|\phi\|_{\cL_{L^\8(\W)}^\8}+\rho_*,\hs &\mb{if } \,r>0;\\[1ex]
e^{-\lam_*t}\|\phi\|_{\cL_{L^\8(\W)}^\8}+\rho_*,\hs &\mb{if }  \,r=0,\\[1ex]
\ea\right.
\ea
\ee
the meaning of  $\phi_k$ is the same as in Lemma \ref{t:1.1}.
\el

\bl\label{t:1.3} Assume that  $\tau$ satisfies {\em(H0)}, $ f$ and $ g$ satisfy {\em (H1)-(H3)}.
Let $q_0<q<\8$, and $h\in L^{\8}(\mathbb{R};L^{q_\a/\a}(\W))\cap C(\mathbb{R};L^{q_\a/\a}(\W))$. Then there exist $M_0,M_1,\lam_1,\rho_1>0$ such that for all $\phi\in \mathscr{X}_1^q$,  we have
\be\label{e:3.6}
\| u_k(t;\phi_k)\|_1^{2} \leq M_0 \|\phi\|_{\cC_{V_1}}^2e^{-\mu_1 t}+M_1\|\phi\|_{\cL_{L^q(\W)}^\8}^{q}e^{-\lam_1 t}+\rho_1^2,\hs \A\,t\geq -r,
\ee
where the meaning of  $\phi_k$ is the same as in Lemma \ref{t:1.1}, $\mu_1$ is the first eigenvalue of the operator $A=-\De$.
   \el
{\bf Proof.}  Since $\{w_j\}^{k}_{j=1}$ be an basis of $W_k:=\mb{span}\{w_1, w_2,\ldots ,w_k\}$, then there exists
$\{d_m^k (t)\}^{k}_{m=1}$ such that $-\De u_k =\sum_{m=1}^{k} d_m^k (t) w_m$. Multiply \eqref{e:2.3} by $d_m^k (t)$ and
sum $m=1,\ldots,k$, we obtain that
\be\label{e:3.23}\begin{split}
\frac{1}{2}\frac{d}{dt}\| u_k\|_1^{2}+\|u_k\|_2^2&=-\( f(u_k)+g(u_k(t-\tau(t,(u_k)_t))+ h_k,\,\De u_k\)\\[1ex]
 &\leq \frac{1}{2}\|u_k\|_2^{2}+ \frac{1}{2}| f(u_k)+g(u_k(t-\tau(t,(u_k)_t))+ h_k|_2^2\\[1ex]
 &\leq \frac{1}{2}\|u_k\|_2^{2}+ C_1\(| f(u_k)|_2^2+|g(u_k(t-\tau(t,(u_k)_t))|_2^2+| h_k|_2^2\).
 \end{split}
 \ee
  By (H1) we deduce that
 \bes
 | f(u_k)|_2^2\leq a_0^2\int_\W(|u_k|^p+1)^2dx\leq  C_2\(|u_k|_q^{2p}+1\),
  \ees
  and
   \bes\begin{split}
  |g(u_k(t-\tau(t,(u_k)_t))|_2^2&\leq b_0^2\int_\W\(\ba{ll} |u_k(t-\tau(t,(u_k)_t)|^{{\b}}+1\ea\)^2dx\\
  &\leq C_3\(\|(u_k)_t\|_{\cL_{L^q(\W)}^\8}^{2\b}+1\).
 \end{split} \ees
Substituting  these estimates into \eqref{e:3.23} we arrive at
\be\label{e:3.24}
\frac{d}{dt}\|u_k\|_1^{2}+ \|u_k\|_2^{2}\leq C_4\(|u_k|_q^{2p}+\|(u_k)_t\|_{\cL_{L^q(\W)}^\8}^{2\b}+1\).
\ee
 Hence by Lemma \ref{t:1.1} one concludes that there exist $M',\lam'>0$ such that
\bes\label{e:3.25}\begin{split}
\frac{d}{dt}\|u_k\|_1^{2}+ \|u_k\|_2^{2}\leq M'\|\phi\|_{\cL_{L^q(\W)}^\8}^qe^{-\lam' t}+C_5,\Hs t>0.
 \end{split}
 \ees
By virtue of the Poinc\'{a}re inequality  we therefore have
\bes\label{e:3.26}\begin{split}
\frac{d}{dt}\|u_k\|_1^{2}&\leq - \|u_k\|_2^{2}+M'\|\phi\|_{\cL_{L^q(\W)}^\8}^qe^{-\lam' t}+C_5\\[1ex]
&\leq - \mu_1\|u_k\|_1^{2}+M'\|\phi\|_{\cL_{L^q(\W)}^\8}^qe^{-\lam_1 t}+C_5,\Hs t>0
 \end{split}
 \ees
where $\lam_1=\min\{\lam',\, \mu_1/2\}$. Thus  by the classical Gronwall lemma one deduces  that  there exist  $M_0, M_1, \rho_1>0$ such   that
\be\label{e:3.27}\begin{split}
\|u_k\|_1^{2}&\leq \|\phi_k\|_{\cC_{ V_1}}^2e^{- \mu_1 t}+\frac{M'}{\mu_1-\lam_1}\|\phi\|_{\cL_{L^q(\W)}^\8}^{q}\(e^{- \lam_1 t}-e^{- \mu_1 t}\)+\rho_1^2\\[1ex]
&\leq M_0\|\phi\|_{\cC_{ V_1}}^2e^{- \mu_1 t}+M_1\|\phi\|_{\cL_{L^q(\W)}^\8}^{q}e^{-\lam_1 t}+\rho_1^2,\Hs t\geq 0.
\end{split}
 \ee
Note that for $t\in [-r,0]$,  the estimate in  \eqref{e:3.6} still holds true because $ u(t)=\phi(t)$ and $e^{-\mu_1 t}\geq 1$. Hence  we  finish  the proof of the lemma. \hfill $\Box$

\br\label{r:1.3}
Let $T>0$. Integrating \eqref{e:3.24} between $0$ and $T$, then by \eqref{e:3.6} one can obtain that
\be\label{e:3.7}
\begin{split}
 \int_0^T\|u_k\|_2^2\,ds&\leq C_T\(\|\phi\|_{\cC_{\bV_1}}^2+\|\phi\|_{\cL_{L^q(\W)}^\8}^{q}+1\).
\end{split}
\ee
\er

\br\label{r:3.1}
If the state-dependent delay in \eqref{e:1.1} is replaced with a constant delay or a time-dependent delay, we can immediately conclude that these estimates in Lemma \ref{t:1.1}-\ref{t:1.3} remain valid for weak solutions by passing to the limit.
However, the delay in \eqref{e:1.1} depends on the state, which poses a more complex situation for us.
\er
\subsection{The proof of Theorems \ref{t:x.1.1} and \ref{t:1.2x}}
We have now gathered enough information to complete the
\Vs
\noindent{{\bf Proof of Theorem \ref{t:x.1.1}.}} For the sake of clarity, we split the argument into several  steps.
\vs
\noindent{\bf Step 1.} Estimates of the nonlinear terms.\vs

 By  the estimates in Lemma \ref{t:1.1} and Lemma \ref{t:1.3}, it can be shown that
\be\label{e:4.1}
 u_k   \in L^{\8}(-r,\8;V_1)\cap L^{\8}(-r,\8;L^q(\W)).
\ee

By  (H1), we have
\bes\begin{split}
\int_{0}^{T}|f(u_k)|_2^2 \,dt&=\int_{0}^{T}\int_{\W}|f(u_k)|^2\,dxdt\\
&\leq \int_{0}^{T}\int_{\W}2 a_0^2 \(|u_k|^{2p}+1\)\,dxdt.\\
\end{split}
\ees
Since $q>q_0\geq \max\{2p,2\b\}$, it can easily be deduced  that there exist $\upsilon_1,\upsilon_2>0$, such that
\be\label{e:3.10}
\int_{\W}2 a_0^2 \(|u_k|^{2p}+1\)\,dx\leq
\upsilon_1 |u_k|_q^q+\upsilon_2|\W|.
\ee
Hence by \eqref{e:4.1}, one can obtain that
$$
f(u_k)\in L^2(0,T;H),\Hs \forall T>0,
$$
and $\{f(u_k)\}$ is uniformly bounded in $L^2(0,T;H)$.
Similarly we have
\be\label{e:3.11}
g(u_k(t-\tau(t,(u_k)_t)) \in L^2(0,T;H),\Hs \forall T>0,
\ee
and $\{g(u_k(t-\tau(t,(u_k)_t))\}$ is uniformly bounded in $L^2(0,T;H)$.

Therefore, one has
\be\label{e:3.12}
w(u_k):=f(u_k)+ g(u_k(t-\tau(t,(u_k)_t))+h_k\in L^2(0,T;H),
\ee
and $\{w(u_k)\}$ is uniformly bounded in $L^2(0,T;H)$.
\vs
\noindent{\bf Step 2.} The convergence of $\{u_k\}$.

 On the one hand, by Lemma \ref{t:1.1} and Lemma \ref{t:1.3}, we can deduce that the  Galerkin approximation  sequence  $\{u_k\}$ is uniformly bounded in $L^{\8}(0,T;V_1)$.

On  the other hand, by \eqref{e:2.2}, we have $u'_k(t)=\sum_{j=1}^k a'_{kj}(t)w_j$, then
multiplying \eqref{e:2.3}  by $a'_{kj}(t)$ and summing over $j$ from $1$ to $k$, one gets that
\be\label{e:3.13}
\begin{split}
|u'_k|_{2}^2+\frac{1}{2}\frac{d}{dt}\int_{\W}\| u_k\|_{1}^2 dx =\( w(u_k),u'_k\)\leq \upsilon_3|w(u_k)|_2^2+ \frac{1}{2} |u'_k|_2^2,
\end{split}\ee
where $\upsilon_3>0$.

Hence we have
\be\label{e:3.14}
|u'_k|_{2}^2+\frac{d}{dt}\int_{\W}\| u_k\|_{1}^2\,dx\leq 2 \upsilon_3|w(u_k)|_2^2.
\ee
Integrate \eqref{e:3.14}  from $0$ to $T$, then
\be\label{e:3.15}
\int_{0}^{T} |u'_k|_{2}^2 \,dt\leq 2\max_{0\leq t\leq T}
\int_{\W}\| u_k(t)\|_{1}^2\,dx+ 2\upsilon_3\int_{0}^{T}|w(u_k)|_2^2  \,dt.
\ee
By \eqref{e:3.6} and \eqref{e:3.12}, one has
\be\label{e:3.16}
u'_k\in L^2(0,T;H),
\ee
and $\{u'_k\}$ is uniformly bounded in $L^2(0,T;H)$.

Thanks to  the Lions-Aubin-Simon Lemma \cite[Corollary 4]{Simon}, it follows from \eqref{e:4.1} and \eqref{e:3.16} that there exists a  subsequence of $\{u_k\}$ converging in $C([0,T];H^{\gamma})$ to a function $u$, where $0 <\gam< \frac{1}{2}$.

Here we assume that $\{u_k\}$ itself is convergent for convenience. i.e., $\{u_k\}$  converges to $u$ in $C([0,T];H^{\gamma})$.
\vs
\noindent{\bf Step 3.} The convergence of $\{f(u_k(x,t))\}$ and $\{g(u_k(t-\tau(t,(u_k)_t))\}$.

 Since $H^{\gamma}\hookrightarrow L^2(\W)$, then in the space $L^2(\W\times(0,T);\mathbb{R})$ one has
\be\label{e:3.17}
u_k(x,t) \rightarrow u(x,t),\Hs a.e.\,(x,t) \in\W\times(0,T).
\ee
Further
\be\label{e:3.18}
f(u_k(x,t)) \rightarrow f(u(x,t)),\Hs a.e.\,(x,t) \in\W\times(0,T),
\ee

On the other hand, by the uniformly boundedness of $\{f(u_k)\}$ in $L^2(0,T;H)$, it is easy to deduce that there exists a weakly convergent subsequence of $\{f(u_k)\}$ that converges to some function $f_0$ (assume that $\{f(u_k)\}$ itself is convergent). By \eqref{e:3.18} and the uniqueness of the limit, we deduce that $f(u(x,t))=f_0$, i.e., $f(u(x,t))$ converges weakly to $f_0$.

We claim that
\be\label{e:4.4}
g(u_k(t-\tau(t,(u_k)_t))  \rightharpoonup g(u(t-\tau(t,u_t))),\Hs k\rightarrow\8.
\ee
Indeed, by \eqref{e:3.6}, one can deduce that there exists a subsequence of $\{u_k\}$ weakly converges to $u$ (assume that $\{u_k\}$ itself is weakly convergent). Then it is easy to see that $\{(u_k)_t\}$ is weakly converges to $u_t$. Hence by (H0), we have
\bes
\tau(t,(u_k)_t)\rightarrow \tau(t, u_t),\Hs k \rightarrow\8.
\ees
Since  $0\leq \tau(\cdot,\cdot)\leq r$, one can obtain that
\bes
u_k(t-\tau(t,(u_k)_t)) \rightarrow u (t-\tau(t,u_t)),\Hs k \rightarrow\8.
\ees

Let
\be\label{e:3.19}
g_k(x,t):=g(u_k(t-\tau(t,(u_k)_t)).
\ee
We can  infer from \eqref{e:3.11} that $g_k$ is bounded in $L^2([0,T];H)$, by which one can deduce that there exists a subsequence of $\{g_k\}$ weakly converges to $g_0$. It may be assumed that $\{g_k\}$   itself converges weakly to $g_0$. Moreover, by $0\leq\tau( \cdot,\cdot)\leq r$ and the continuity of $g$, one has
\be\label{e:3.20}
g(u_k(x,t-\tau(t,(u_k)_t)))\rightarrow g(u(x,t-\tau(t,u_t))),\Hs a.e.\,(x,t) \in\W\times(0,T).
\ee

Therefore by \eqref{e:3.19} and \eqref{e:3.20}, it can be seen that
\bes\label{e:3.21}
g_0(x,t)=g(u (x,t-\tau(t,u_t))).
\ees
That's to say, $g(u_k(t-\tau(t,(u_k)_t)))$ weakly converges to  $g(u(t-\tau(t,u_t)))$ as $k\ra\8$.
\vs
\noindent{\bf Step 4.}  Verification of the conclusion.

 Furthermore, it is easy to see that $h_k$ converges to $h$ as $k \ra\8$.  Passing to the limit one immediately arrive at
 \be\label{e:4.9}
\left\{
\begin{split}&\(\frac{d u }{dt}-\De u ,w_j\)=\( f( u )+g(( u (t-\tau(t,u_t)))+ h ,w_j\),\Hs j=1,2,\ldots,\8\\
 & u (s)=\phi (s),\hs s\in[-r,0].
\end{split}
\right.
\ee
Then   repeating the discussion in Lemma \ref{t:1.1} and Lemma \ref{t:1.3}, there exist positive constants $B_i$, $\eta_i$ and $\rho_i$ ($i=0,1$) such that \eqref{e:q.3} and \eqref{e:q.5} hold true.

Similar to \eqref{e:3.1} and \eqref{e:3.6}, it follows that $u\in L^\8\(-r,\8; V_1\)\cap L^\8\(-r,\8;L^q(\W)\)$. Repeat from \eqref{e:4.1} to \eqref{e:3.20}, it can be deduced that $u'\in L^2(0,T;H)$. Note that by replacing $u_k$ with $u$, \eqref{e:3.7} still holds, so $u\in L^2(0,T;V_2)$.

The continuity of $u$ in $V_1$ now follows from  \cite[Chap. II, Theorem 3.3]{Tem} on abstract linear equations. This completes  the proof of theorem.   \hfill $\Box$

\Vs
\noindent{{\bf Proof of Theorem \ref{t:1.2x}.}}
 The conclusions about existence and regularity clearly hold. Passing to limits as $k\rightarrow\8$, then
the estimates in \eqref{e:q.4} is the result of   Lemma \ref{t:1.2} and Theorem \ref{t:x.1.1}.
\hfill$\Box$

\section{Fundamental properties in fractional power spaces}\label{s:4}
Assume  $\tau$ satisfies (H0),  $ f$ and $ g$ satisfy (H1)-(H3), and let $q_0<q\leq \8$.  Given $\phi\in\sX_1^q$ and $ h\in L^{\8}(\mathbb{R};L^{q_\a/\a}(\W))\cap C(\mathbb{R};L^{q_\a/\a}(\W))$, denote by $\sS(\phi)$ the set of solutions of \eqref{e:xx2} with initial value $\phi$. Define the process $\Phi$ from $\R^+\X \R \X\sX_1^q$ to $\sX_1^q $ as below:
\be\label{e:4.13}
\Phi=\Phi(t,\sigma,\phi)=\{\bu_t\,|\,\bu\in \sS(\phi)\}.
\ee
By Remark \ref{r:1.1} it is trivial to see that  $\Phi$ has the semigroup property:
\benu
\item[(i)] $\Phi(0,\sigma,\phi)=\phi$ for all $\sigma\in\R$, $\phi\in\sX_1^q$; and\vs
\item[(ii)] $\Phi(t+s,\sigma,\phi)=\Phi\(t,s+\sigma,\Phi(s,\sigma,\phi)\)$ for all $t,s\geq 0$ and $\phi\in\sX_1^q$.
\eenu

In the Section \ref{s:3}, we obtain some estimates of the global dissipativity of the solution of \eqref{e:1.1}, which remain true for \eqref{e:xx2} for all $\sigma\in\R$, with the constants therein being independent of $\sigma$.
Now we give results on boundedness and equi-continuity in fractional power spaces, which are the necessary for us to consider the pullback asymptotic compactness of $\Phi$.
For convenience, we rewrite \eqref{e:xx2} as follows:
\be\label{e:5.1}
u'+A u=w.
\ee
where 
 $w(t)=f(u(t))+g(u(t-\tau(t+\sigma,u_t)))+h(t+\sigma, x)$. By Theorem \ref{t:x.1.1} and (H1), one has $w\in L^{\8}(0,\8;H)$. Hence we assume that there exists a positive  constant $b$, such that
\be\label{e:5.2}
\|w\|_{L^{\8}(0,\8;H)}\leq b.
\ee

In this section, we renormalize $\|\cdot\|$ to denote the norm in $H=L^2(\W)$. Let $\|\cdot\|_{\zeta}$ denote the norm in the fractional power space $H^{\zeta}$ if $0<\zeta<1$, and it is worth noting that $H^{\frac{1}{2}}=V_1=H_0^1(\W)$, i.e., $\|\cdot\| _{V_1}=\|\cdot\|_{\frac{1}{2}}$.

\bt\label{l:1.1}Suppose $\tau$, $f$ and $g$ satisfy {\em(H0)}-{\em(H3)}. Let
$\frac{1}{2}< \zeta <1$, $q_0<q<\8$ and $ h\in L^{\8}(\mathbb{R};L^{q_\a/\a}(\W))\cap C(\mathbb{R};L^{q_\a/\a}(\W))$. Then for any $\sigma\in\R$, $\phi\in \mathscr{X}_1^q$,  there exist  $M^*=M^*(\zeta),C_b=C(b,\zeta)>0$ such that for any solution $u$ of \eqref{e:5.1},
we have
\be\label{e:5.3}
\|u\|_{L^{\8}(r,\8;H^{\zeta})}\leq M^* \|\phi\|_{\cL_{L^q(\W)}^\8} + C_b.
\ee
\et
{\bf Proof.} By Theorem \ref{t:2.2}, we deduce that
\be\label{e:5.4}
u(t)=T(t)u(0)+\int_0^t T(t-s)w(s)ds, \Hs t\geq r.
\ee
where $T(t)=e^{-tA}$ is the $C_0$ semigroup generated by $A$. Then we have
\be\label{e:5.5}
\begin{split}
\|u(t)\|_{\zeta} &=\left\|A^{\zeta}u(t)\right\|=\left\|A^{\zeta}T(t)u(0)+\int_0^t A^{\zeta} T(t-s)w(s)ds\right\|\\
&\leq \left\|A^{\zeta}T(t)u(0)\right\|+\int_0^t\left\| A^{\zeta} T(t-s)w(s)\right\|ds\\
&=I_0+I_1,
\end{split}\ee
here we use the fact that $A$ is a densely defined closed operator.

By Lemma \ref{l:1.2} one has
\be\label{e:5.6}
\begin{split}
I_0=&\left\|A^{\zeta}T(t)u(0)\right\|=\left\|A^{\zeta}T(t)\phi(0)\right\|\\& \leq C_{\zeta}r^{-\zeta} e^{-r \zeta}\|u_t\|_{L^{\8}(-r,0;H)}\\
&=C_{\zeta,r}\|u_t\|_{L^{\8}(-r,0;H)}
,\Hs t\geq r,
\end{split}
\ee
where $C_{ \zeta,r} =C_{\zeta}r^{-\zeta} e^{-r \zeta}$ is a constant relies on $\zeta$ and $r$.

On the other hand, by \eqref{e:1.14} we observe that
\be\label{e:5.7}
\begin{split}
I_1=&\int_0^t \left\|A^{\zeta} T(t-s)w(s)\right\|ds\\
&\leq C_{\zeta} \int_0^t(t-s)^{-\zeta}e^{-\delta(t-s)}\|w(s)\|_{L^{\8}(0,\8;H)} ds\\
& \leq(\mb{by \eqref{e:5.2}}) \leq b\, C_{\zeta} \int_0^t  (t-s)^{-\zeta}e^{-\delta(t-s)}ds.
\end{split}\ee
We need to evaluate the integral $\int_0^t  (t-s)^{-\zeta}e^{-\delta(t-s)}ds$:
\bes\label{e:5.8}
\begin{split}
\int_0^t  (t-s)^{-\zeta}e^{-\delta(t-s)}ds &\xlongequal{t-s=\eta}\int_0^t \eta^{-\zeta}e^{-\delta \eta}d\eta\xlongequal{\delta\eta=\gamma} \int_0^{\delta t}(\frac{\gamma}{\delta})^{-\zeta}e^{-\gamma}\frac{1}{\delta}d\gamma\\
&=(\frac{1}{\delta})^{1-\zeta}\int_0^{\delta t}\gamma^{-\zeta}e^{-\gamma}d\gamma\leq \delta^{\zeta-1}\int_0^{\8}\gamma^{-\zeta}e^{-\gamma}d\gamma\\
&=\delta^{\zeta-1} \Gamma(1-\zeta),
\end{split}
\ees
where $\Gamma(\cdot)$ is the Gamma Function.

Combining this with \eqref{e:5.7}, we have
\be\label{e:5.9}
\int_0^t \left\|A^{\zeta} T(t-s)w(s)\right\|ds\leq  b \,C_{\zeta}\delta^{\zeta-1}\Gamma(1-\zeta).
\ee
Substituting \eqref{e:5.6} and \eqref{e:5.9} into \eqref{e:5.5}, one can deduce that
\be\label{e:5.10}
\|u(t)\|_{\zeta}\leq C_{\zeta,r}\|u_t\|_{L^{\8}(-r,0;H)}+b \,C_{\zeta}\delta^{\zeta-1}\Gamma(1-\zeta),\Hs t\geq r.
\ee
Hence by \eqref{e:5.10} and Theorem \ref{t:x.1.1}, we finally arrive at the estimate in \eqref{e:5.3}. \hfill $\Box$

Let
\be\label{e:5.21}
\mathscr{B}^q(R)=\{\phi \in \mathscr{X}^q_1\,|\, \|\phi\|_{\mathcal{C}_{V_1}}+\|\phi\|_{\mathcal{L}^{\8}_{L^q (\W)}}\leq R\}.
\ee
We can prove an equi-continuity result concerning solutions of  \eqref{e:5.1}.
\bt\label{l:2.1}
Assume the hypotheses in Theorem \ref{l:1.1}. Let $\frac{1}{2}<\zeta<1$ and $B\subset \sB^q(R)$, then for any $\eta>0$, there exists $L=L(\eta)>0$ such that for all $u\in  \sS (B):= \Cup_{\phi\in B}\sS(\phi)$, we have
\bes\label{e:5.11}
\|u(t+\nu)-u(t)\|_{\frac{1}{2}}\leq L \, (\nu^{\zeta-\frac{1}{2}}+\nu^{\frac{1}{2}}),\Hs \forall \nu>0,\hs t\geq 2r+\eta.
\ees
\et
{\bf Proof.} Consider
\be\label{e:5.12}
u(t)=T(t-2r)u(2r)+\int^t_{2r}T(t-s)w(s)ds,\Hs t\geq 2r+\eta.
\ee
Then we have
\be\label{e:5.13}
\begin{split}
\left\|u(t+\nu)-u(t)\right\|_{\frac{1}{2}}&\leq\|T(t-2r)(T(\nu)-I)u(2r)\|_{\frac{1}{2}}\\
&\hs+\left\|\int^t_{2r}T(t-s)w(s)ds-\int^{t+\nu}_{2r}T(t+\nu-s)w(s)ds\right\|_{\frac{1}{2}}.
\end{split}\ee

First we evaluate $\|T(t-2r) (T(\tau)-I) u(2r)\|_{\frac{1}{2}}$. It is easy to see that
$$
\|T(t)\|\leq M_3(s),\Hs \forall t\geq s,
$$
then we have
\be\label{e:5.19}
\begin{split}
\|T(t-2r)(T(\nu)-I)u(2r)\|_{\frac{1}{2}}&= \|A^{\frac{1}{2}}T(t-2r)(T(\nu)-I)u(2r) \|\\
&\leq M_3(\eta)\,\|A^{\frac{1}{2}}(T(\nu)-I)u(2r)\|\\
&\leq M_3(\eta)\,\|A^{-(\zeta-\frac{1}{2})}(T(\nu)-I) A^{\zeta}u(2r)\|,
\end{split}\ee
where $\frac{1}{2}<\zeta<1$.
We also observe that
\be\label{e:5.16}
\|A^{-(\zeta-\frac{1}{2}) }(e^{-tA} -1)\|\leq \left\|\int_0^t A^{\frac{3}{2}-\zeta}e^{-sA}ds\right\|\leq C_{\frac{3}{2}-\zeta}\int^t_0 s^{\zeta-\frac{3}{2}}ds\leq  C'_{\frac{3}{2}-\zeta} t^{\zeta-\frac{1}{2}},\ee
where the constants $C_{\frac{3}{2}-\zeta}$ and $C'_{\frac{3}{2}-\zeta}$ are both related to $\frac{3}{2}-\zeta$. Substituting this into \eqref{e:5.19}, one can obtain that
\be\label{e:5.14}
\|T(t-2r)(T(\nu)-I)u(2r)\|_{\frac{1}{2}}\leq M_4(\eta)\,  \nu^{\zeta-\frac{1}{2}}\|u(2r)\|_{\zeta},\Hs t\geq 2r+\eta.
\ee

Let us  evaluate the second term on the right side of \eqref{e:5.13}, obviously 
\be\label{e:5.15}
\left\|\int^t_{2r}T(t-s)w(s)ds-\int^{t+\nu}_{2r}T(t+\nu-s)w(s)ds\right\|_{\frac{1}{2}}\leq  I_0+I_1,
\ee
where
\bes
\begin{split}
I_0&=\left\|\int^t_{2r}(T(t-s)-T(t+\nu-s))w(s)ds\right\|_{\frac{1}{2}},\\
I_1&=\left\|\int^{t+\nu}_{t}T(t+\nu-s)w(s)ds\right\|_{\frac{1}{2}}.
\end{split}\ees

By \eqref{e:5.16} and \eqref{e:1.14}, we deduce that
\bes\label{e:5.17}
\begin{split}
I_0&=\left\|\int^t_{2r}A^{\frac{1}{2}} T(t-s)(I-T(\nu)) w(s)ds\right\|\\
&=\left\|\int^t_{2r}A^{-(\zeta-\frac{1}{2})} (I-T(\nu))A^{\zeta}T(t-s) w(s)ds\right\|\\
&\leq\int^t_{2r} \left\|A^{-(\zeta-\frac{1}{2})}(I-T(\nu))\right\|\cdot \left\|A^{\zeta}T(t-s)\right\| \cdot \left\|w(s)\right\|ds\\
&\leq C'_{\frac{3}{2}-\zeta}\nu^{\zeta-\frac{1}{2}}C_{\zeta} \|w(s)\|_{L^{\8}(0,\8;H)}\int^t_{2r}(t-s)^{-\zeta}e^{-\delta(t-s)}ds\\
&\leq C'_{\frac{3}{2}-\zeta}\nu^{\zeta-\frac{1}{2}}C_{\zeta} b\, \delta^{\zeta-1}\Gamma(1-\zeta),
\end{split}
\ees
and
\bes\label{e:5.18}
\begin{split}
I_1&=\left\|\int^{t+\nu}_{t}A^{\frac{1}{2}}T(t+\nu-s)w(s)ds\right\|\\
&\leq \int^{t+\nu}_{t} C_{\frac{1}{2}}(t+\nu-s)^{-\frac{1}{2}} e^{-(t+\nu-s)\delta}\|w(s)\|_{L^{\8}(0,\8;H)}ds\\
&\leq  C_{\frac{1}{2}}\|w(s)\|_{L^{\8}(0,\8;H)}\int^{t+\nu}_{t}(t+\nu-s)^{-\frac{1}{2}}ds\\
&\leq 2\,C_{\frac{1}{2}} b\,\nu^{\frac{1}{2}}.
\end{split}\ees
Substituting these estimates into \eqref{e:5.15}, one conclude that
\be\label{e:5.20}
\begin{split}
&\left\|\int^t_{2r}T(t-s)w(s)ds-\int^{t+\nu}_{2r}T(t+\nu-s)w(s)ds\right\|_{\frac{1}{2}}\leq  M_5 (\nu^{\zeta-\frac{1}{2}}+\nu^{\frac{1}{2}}).
\end{split}
\ee
Thus the conclusion can be obtained naturally  by Theorem \ref{l:1.1},  \eqref{e:5.14}  and \eqref{e:5.20}.\hfill $\Box$

As a corollary to the above two theorems, we have

\bl\label{l:2.2}
Assume the hypotheses in Theorem \ref{l:1.1}.
Let  $B\subset\sB^q(R)$, and  $\{u^k\}\subset \sS(B)$ be a sequence of solutions of \eqref{e:xx2}. Then for any compact interval $[\eta,T]\subset (2r,\8)$, there is a subsequence of $\{u^k\}$ converging  in $C([\eta,T]; V_1)$ to a solution $u$ of \eqref{e:xx2}.
\el
{\bf Proof.} Theorem \ref{l:2.1} shows that the sequence   $\{u^k\}$ is equi-continuous on $[\eta,T]$ for all $\sigma\in\R$. By Theorem \ref{l:1.1}, we see that there is a bounded set $\tilde{B}$ in $H^{\zeta}$($\frac{1}{2}<\zeta<1$) such that
$$
 u^k(t)\in \tilde{B},\Hs   t\in [\eta,T],\,\,k=1,2,\cdots.
$$
Since $ H^{\zeta}\hookrightarrow\hookrightarrow V_1$, we can obtain that $\tilde{B}$ is precompact in $ V_1$. By virtue of the Arzela-Ascoli's theorem, there is a subsequence $\{ u^{k_i}\}$ converging   in $C([\eta,T]; V_1)$ to a function $ u$.  On the other hand, using the estimates given in previous sections it can be shown by standard argument that, up to a subsequence, $\{ u^{k_i}\}$ converges in appropriate spaces to a (weak) solution $\~u$ of equation \eqref{e:xx2} on $[\eta,T]$. As $ u^{k_i}\ra u$ ($k\ra\8$) in $C([\eta,T]; V_1)$, one naturally has $\~u=u$. The proof is complete.  \hfill$\Box$

\section{Global pullback $\sD$-attractor of system \eqref{e:1.1}}\label{s:5}
It has already been mentioned that the uniqueness of the solution in \eqref{e:1.1}  is hard to obtain, i.e. the process $\Phi(t,\sigma,\phi)$ is actually set-valued.
This section is concerned with the existence of the global pullback $\sD$-attractor of the set-valued process $\Phi(t,\sigma,\phi)$.

\subsection{Pullback $\sD$-attractors of set-valued processes}
Let $(X,d)$ be a metric space, $P(X)$ the family of all nonempty subsets of $X$, and $d_H(A,B)$ the {\em Hausdorff  semi-distance} between $A$ and $B$, i.e.,
$$
d_H(A,B)=\sup_{x\in A}\inf_{y\in B}d(x,y),\Hs \forall A,B\in P(X).
$$
\begin{definition}\label{d:5.1x}
A mapping $\Phi:\R^+\X \R \X X \rightarrow P(X)$ is called a set-valued process if
\benu
\item[(i)] $\Phi(0,\sigma,\phi)=\{\phi\}$ for all $\sigma\in\R$, $\phi\in X$; and\vs
\item[(ii)] $\Phi(t+s,\sigma,\phi)\subset \Phi\(t,s+\sigma,\Phi(s,\sigma,\phi)\)$ for all $t,s\geq 0$ and $\phi\in X$, where
    $$
    \Phi(t,s,A)=\bigcup_{y\in A} \Phi(t,s,y),\Hs \forall A\in P(X).
      $$
\eenu
In particular, $\Phi$ is called strict if $\Phi(t+s,\sigma,\phi)= \Phi\(t,s,\Phi(s,\sigma,\phi)\)$.
\end{definition}

\begin{remark}
By Remark \ref{r:1.1} it is easy to see that the set-valued process associate with the system \eqref{e:1.1} is strict, and all of our work below is based on this premise.
\end{remark}

\begin{definition}\label{d:5.2x}
A family of nonempty sets $\cD:=\{\cD(\sigma)\subset X\}_{\sigma\in\R}$ will be called a nonautonomous sets.
\end{definition}

Let $\Phi$ be the set-valued process on $X$ and $\sD$ a  family of nonautonomous sets. Then we borrow some ideas from \cite{Car1} and \cite{Coti} to establish the framework of pullback attractor in the context
of set-valued dynamical system. For convenience we rewrite $\Phi(t,\sigma,\phi)=\Phi(t,\sigma)\phi$.

\begin{definition}\label{d:5.3}
A set-valued process $\Phi$ is called pullback $\sD$-asymptotically compact, if for any $\cD\in\sD$, $\sigma\in\R$, and any sequence $t_n\ra +\8$ and $\phi_n\in\cD(\sigma-t_n)$, the sequence $\Phi(t_n,\sigma-t_n)\phi_n$ has a convergent subsequence.
\end{definition}

\begin{definition}
A set-valued process $\Phi$ is called pullback $\sD$-dissipative if there exists a nonautonomous set $\cB$ such that for any $\cD\in\sD$ and $\sigma\in\R$, there exists $t_0=t_0(\cD,\sigma)$ such that
$$
\Phi(t,\sigma-t)\cD(\sigma-t)\subset \cB(\sigma), \Hs \forall t\geq t_0.
$$
The nonautonomous set $\cB$ is called pullback $\sD$-absorbing.
\end{definition}

\begin{definition}\label{d:5.4}
A  nonautonomous set $\cA$ is called a global pullback $\sD$-attractor for the set-valued process if it fulfills the following requirements:
\benu
\item[(1)] $\cA(\sigma)$ is compact for each $\sigma\in \R$;\vs
\item[(2)] $\cA(\sigma)$ is pullback $\sD$-attracting for each $\sigma\in \R$, i.e., for any $\cD\in\sD$,
    $$
    \lim_{t\ra+\8}d_H\(\Phi(t,\sig-t)\cD(\sig-t),\, \cA(\sig)\)=0,\Hs \sig\in\R.
    $$
\item[(3)] if the nonautonomous set $\cF$ is pullback $\sD$-attracting with $\cF(\sigma)$ being closed, then $\cA(\sigma) \subset \cF(\sigma)$ for each $\sigma\in\R$.
\eenu
\end{definition}

\begin{theorem}\label{t:5.1x}
Suppose the set-valued process $\Phi$ is pullback $\sD$-asymptotically compact, and there exists a set $\cB\in\sD$ which is pullback $\sD$-absorbing. Then $\Phi$ has a global unique pullback $\sD$-attractor $\cA$ given by
\be\label{e:5.1*}
\cA(\sig)=\Cap_{\tau\geq 0}\ol{\Cup_{t\geq \tau}\Phi(t,\sig-t)\cB(\sig-t)},\Hs \sig\in\R.
\ee
\end{theorem}
{\bf Proof.} Let $\cB=\{\cB(\sigma)\,|\,\sigma\in\R\}$ be a monotonic pullback $\sD$-absorbing set. By \cite[pp. 489, Proposition 12]{Car1} it is can be easily deduce that $$\cA(\sigma):=\Cap_{\tau\geq 0}\ol{\Cup_{t\geq \tau}\Phi(t,\sig-t)\cB(\sig-t)}
$$
is nonempty and compact and
\be\label{e:5.2*}
\lim_{t\ra+\8}d_H\(\Phi(t,\sig-t)\cD(\sig-t),\, \cA(\sig)\)=0,\Hs \sig\in\R.
\ee

We need to show that $\cA=\{\cA(\sigma)\,|\,\sigma\in\R\}$ is pullback $\sD$-attracting. For any $s\in\R^{+},\sigma\in\R$ with $s>\sigma$,
we can find $ \bar{\nu}>s$ such that
$$
\Phi(\nu,\sigma-\nu)\cD(\sigma-\nu)\subset
 \Phi(s,\sigma-s,\Phi(\nu-s,\sigma-\nu)\cD(\sigma-\nu))\subset
 \Phi(s,\sigma-s)\cD(\sigma-s),\,\,\,\, \forall \nu\geq\bar{\nu}.
 $$
Therefore, we have
$$
d_H\(\Phi(\nu,\sigma-\nu)\cD(\sigma-\nu),\cA(\sigma) \)\leq d_H\(\Phi(s,\sigma-s)\cD(\sigma-s),\cA(\sigma)\).
$$
Further, we can obtain that
$$
\lim_{\nu\ra+\8}\sup d_H\(\Phi(\nu,\sigma-\nu)\cD(\sigma-\nu),\cA(\sigma) \)\leq d_H\(\Phi(s,\sigma-s)\cD(\sigma-s),\cA(\sigma)\),
$$
Then in light of \eqref{e:5.2*}, one immediately concludes the pullback attraction property by taking the limit $s\ra +\8$ in the above inequality.

For minimality, let $\cF=\{\cF(\sigma)\,|\,\sigma\in\R\}$ be another closed pullback $\sD$-attracting set, i.e.,
$$
\lim_{t\ra+\8}d_H\(\Phi(t,\sig-t)\cD(\sig-t),\, \cF(\sig)\)=0.
$$
Let $\eta\in\cA(\sigma)$, then by some basic knowledge in \cite{Car1} and \cite{Coti} we have $\eta=\lim_{n\ra\8} \Phi(t_n,\sigma-t_n) y_n$ for  some sequences $t_n \ra+\8$ and $y_n\in\cD(\sigma-t_n)$. Consequently one has $\eta\in\cF(\sigma)$. Thus,  $\cA(\sigma)\in \cF(\sigma)$.

The uniqueness of $\sD^q$-attractor of $\Phi$ follows from the minimality requirement in the definition of $\sD^q$-attractors. The proof of the theorem is complete.
\hfill$\Box$

\subsection{Global pullback $\sD$-attractor of \eqref{e:1.1}}

It is evident that the system \eqref{e:1.1} is actually nonautonomous, and in order to describe the global dynamical behaviour of such a system, we study it in a pullback framework. Therefore we need to consider the corresponding system given in \eqref{e:xx2} with the time symbol $\sigma\in\R$. As we mentioned earlier, Theorems \ref{t:x.1.1} and \ref{t:1.2x} as well as all estimates given in Sections \ref{s:3}-\ref{s:4} hold for \eqref{e:xx2} for all $\sigma\in\R$, since all the constants involved in estimates are independent of $\sigma$.

let $(X,d)$ be the  space $(\sX_1^q,d)$, where the metric $\mb{d}(\cdot,\cdot)$ is induced by the norm of $\cC_{ V_1}$:
$$
d(\phi,\psi)=\|\phi-\psi\|_{\cC_{ V_1}}:=\max_{s\in[-r,0]}\|\phi-\psi\|_{1},\Hs \forall \phi,\psi\in\sX_1^q.
$$
Recall that we define the (strict) set-valued process $\Phi$ on $\sX_1^q$ in \eqref{e:4.13}:
$$
\Phi(t,\sigma,\phi)=\Phi(t,\sigma)\phi=\{\bu_t\,|\,\bu\in \sS(\phi)\},
\Hs t\geq 0,\,\sigma\in\R,
$$
where $\sS(\phi)$ is the set of solutions of \eqref{e:xx2} with initial value $\phi$.

Here we identify a nonempty set $D\subset \sX_1^q$ and the nonautonomous one $\cD:=\{\cD(\sigma)\}_{\sigma\in\R}$ in $\sX_1^q$ with $\cD(\sigma)\equiv D$. Define
$$
\sD:=\{\sB^q(R)\,|\,R>0\},
$$
where $\sB^q(R)$ can be find in \eqref{e:5.21}.

\bt\label{t:5.1}
Assume  $\tau$ satisfies (H0),  $ f$ and $ g$ satisfy (H1)-(H4). Let $q_0<q\leq \8$, and $ h\in L^{\8}(\mathbb{R};L^{q_\a/\a}(\W))\cap C(\mathbb{R};L^{q_\a/\a}(\W))$. Then
$\Phi$ has a global unique pullback $\sD$-attractor $\cA=\{\cA(\sigma)\}_{\sigma\in\R}$ in $\sX_1^q$ which possesses some properties as stated in Definition \ref{d:5.4}.
\et
 {\bf Proof.}
It is trivial to check that all the aforementioned estimates  remain true for the system \eqref{e:xx2}.
Let $R_0=\max\{\rho_0,\rho_1\}+1$, where $\rho_0$ and $\rho_1$ are  positive numbers respectively given in Theorems \ref{t:x.1.1}. Then one easily deduces that $\cB:=\sB^q(R_0)$ is a pullback $\sD$-absorbing set.
To apply Theorem \ref{t:5.1x} to establish the existence of a global pullback attractor, it is necessary to verify the pullback asymptotic compactness of $\Phi$.

Let $R>0$ and $\sigma\in\R$ be given arbitrary. Let $\phi_n$ be a sequence in $\sB^q(R)$ and $t_n\ra+\8$.
Denote $\{u^k_n\}_k$ be the sequence generated by the set-valued process $\Phi(t_n,\sigma-t_n)\phi_n=\{u(t_n+s)\,|\,u\in \sS(\phi_n)\}$, where $s\in[-r,0]$. By two estimates in Theorem \ref{t:x.1.1} there exists $t_0>0$ such that
\be\label{e:5.22}
\Phi(t_n,\sigma-t_n)\phi_n\subset \sB^q(R), \Hs \forall t_n\geq t_0.
\ee

On the other hand, for any $i\geq 2$, we have
\be\label{e:5.23}
\Phi(t_i,\sigma-t_i)\phi_i=\Phi(t_1,\sigma-t_1,\Phi(t_i-t_1,\sigma-t_i,\phi_i))
\ee
Since $t_n\ra+\8$, we can assume that there exists $i_0\geq 2$ such that
$$
t_i-t_1 \geq t_0,\Hs i\geq i_0.
$$
Combining this with \eqref{e:5.22} and \eqref{e:5.23}, then
Lemma \ref{l:2.2} implies that the sequence $\{u^k_n\}_k$ has a subsequence $\{u^{k_j}_{n_i}\}_{k_j}$ which is convergent in $\cC_{\bV_1}$.
Further by Theorem \ref{t:5.1x}, we finally arrive at the conclusion in Theorem \ref{t:5.1}.
\hfill$\Box$

\bt\label{t:5.2} Assume hypotheses  {\em (H0)-(H4)} are fulfilled. Let  $\bh\in  L^{\8}(\mathbb{R};L^{\8}(\W))\cap C(\mathbb{R};L^{\8}(\W))$,
then  $\Phi$  has a  unique global pullback $\sD^\8$-attractor $\tilde{\cA}$ in $\sX_1^{\8}$.
\et

\section*{References}
\Vs
\end{document}